\documentclass{amsart2000} 
\usepackage{graphicx,pstricks, amssymb,  amscd2000}
\swapnumbers \newtheorem*{theorem*}{Theorem}
\newtheorem{theorem}{Theorem}[section]
\newtheorem{lemma}[theorem]{Lemma}

\newtheorem{corollary}[theorem]{Corollary}

\newtheorem{assumption}[theorem]{Assumption}

\theoremstyle{definition} 
\newtheorem{definition}[theorem]{Definition}
\newtheorem{example}[theorem]{Example}

\theoremstyle{remark} 
\newtheorem{remark}[theorem]{Remark}

\numberwithin{equation}{section}
\newcommand{\firef}[1]{Figure~{\rm\ref{#1}}}
\newcommand{\thref}[1]{Theorem~{\rm\ref{#1}}}
\newcommand{\asref}[1]{Assumption~{\rm\ref{#1}}}

\newcommand{\leref}[1]{Lemma~{\rm\ref{#1}}}
\newcommand{\coref}[1]{Corollary~{\rm\ref{#1}}}

\newcommand{\exref}[1]{Example~{\rm\ref{#1}}}

\newcommand{\seref}[1]{Section~{\rm\ref{#1}}}

\newcommand{\fig}[1]
{\raisebox{-0.5\height}%
  {\includegraphics{#1}}}

\newcommand{\newfig}[2]
{\raisebox{-0.5\height}%
{\includegraphics{#1}{#2}}
}

\newcommand{\st}{\; | \;}                     
\newcommand{\ttt}{\otimes}                    
\newcommand{\sttt}{\!\otimes\!}               
\newcommand{\tta}{\otimes_A}                  
\newcommand{\tbox}{\boxtimes} 

\newcommand{\<}{\langle} 
\renewcommand{\>}{\rangle}
\newcommand{\surjto}{\twoheadrightarrow}      
\newcommand{\injto}{\hookrightarrow}          
\newcommand{\isoto}{\xrightarrow{\sim}}       
\newcommand{\xxto}{\xrightarrow}              
\newcommand{\wt}{\mathrm{wt}}                 
\newcommand{\twt}{\widetilde{\mathrm{wt}}}    

\newcommand{\one}{\mathbf{1}}     
\newcommand{\Cset}{\mathbb{C}}    
\newcommand{\Ctimes}{\mathbb{C}^\times}       
\newcommand{\V}{{\mathcal{V}}}    
\newcommand{\C}{\mathcal{C}}      
\newcommand{\A}{\mathcal{A}}      
\newcommand{\F}{\mathcal{F}}      
\renewcommand{\vec}{\mathcal{V}ec}
\newcommand{\tA}{\tilde{A}} 
\newcommand{\mut}{\widetilde{\mu}}
\newcommand{\cd}{{\C[D(G)]}}      
\newcommand{\rcd}{\check R^D}

\newcommand{\ga}{\gamma} 

\newcommand{\de}{\delta} 
\newcommand{\De}{\Delta}
\newcommand{\la}{\lambda} 
\newcommand{\ph}{\varphi}
\newcommand{\Ph}{\Phi} 

\renewcommand{\th}{\theta} 
\newcommand{\om}{\omega}
\newcommand{\Ghat}{\widehat{G}}       
\newcommand{\tdg}{\widehat{D(G)}}     
\newcommand{\slthat}{\widehat{\mathfrak{sl}}_2}    

\DeclareMathOperator{\Rep}{Rep}
 
\DeclareMathOperator{\id}{id}
 
\DeclareMathOperator{\Hom}{Hom}
 
\DeclareMathOperator{\End}{End}
\DeclareMathOperator{\tr}{tr}

\begin{document}
\title{Modular  categories and orbifold models} 
\author{Alexander Kirillov, Jr.}
\address{Department of Mathematics, SUNY at Stony Brook, 
         Stony Brook,  NY 11794, USA} 
\email{kirillov@math.sunysb.edu}
\urladdr{http://www.math.sunysb.edu/\textasciitilde kirillov/}
\thanks{The author was supported in part by NSF grant DMS9970473.}
                           
\begin{abstract}
  In this paper, we try to answer the following question: given a
  modular tensor category $\A$ with an action of a compact group $G$,
  is it possible to describe in a suitable sense the ``quotient''
  category $\A/G$? We give a full answer in the case when $\A=\vec$ is
  the category of vector spaces; in this case, $\vec/G$ turns out to
  be the category of representation of Drinfeld's double $D(G)$. This
  should be considered as category theory analog of topological
  identity $\{pt\}/\!/G=BG$.  

  This implies  a conjecture of Dijkgraaf, Vafa, E.~Verlinde
  and H.~Verlinde regarding so-called orbifold conformal field theories:
  if $\V$ is a vertex operator algebra which has a unique irreducible
  module, $\V$ itself, and $G$ is a compact group of automorphisms of
  $\V$, and some not too restrictive technical conditions are
  satisfied, then $G$ is finite, and the category of representations
  of the algebra of invariants, $\V^G$, is equivalent as a tensor
  category to the category of representations of Drinfeld's double
  $D(G)$. We also get some partial results in the non-holomorphic
  case, i.e. when $\V$ has more than one simple module. 
\end{abstract}

\maketitle
\section*{Introduction}
The goal of this paper is to discuss the properties of the so-called
orbifold models of Conformal Field Theory from the categorical point
of view.

For readers convenience, we recall here the main definitions and
results, assuming that the reader is familiar with the notion of a
vertex operator algebra.  Let $\V$ be a VOA and $G$ --- a finite group
acting on $\V$ by automorphisms. Then the subspace of invariants $\V^G$
is itself a VOA. The main question is: is it possible to describe the
category of $\V^G$-modules in terms of the category of $\V$-modules
and the group $G$?

This question was asked in this form in \cite{DVVV}. They didn't give
a full answer, but did suggest (without a proof) an answer in a
special case, when $\V$ is holomorphic, i.e. has only one simple
module (vacuum module). In this case, the category of representations
of $\V$ is equivalent to the category of vector spaces. The answer
suggested in \cite{DVVV} and further discussed in \cite{DPR} is that
in this case, the category of $\V^G$-modules is equivalent to the
category of modules over the (twisted) Drinfeld double $D(G)$ of the
group $G$. The case of VOA's coming from Wess-Zumino-Witten model (or,
equivalently, from affine Lie algebras) was studied in detail in
\cite{KT}.

Many of the results of \cite{DVVV} were rigorously proved in the
language of VOA's in a series of papers of Dong, Li, and Mason; in
particular, it is proved in \cite{DLM1} that in the holomorphic case,
$\V$ considered as a module over both $\V^G$ and $G$ can be written as
\begin{equation}\label{dlm1}
\V=\bigoplus_{\la\in \Ghat} V_\la \ttt M_\la
\end{equation}
where $V_\la$ are irreducible $G$-modules and $M_\la$ are non-zero
simple pairwise non-isomorphic $\V^G$-modules. However, even in the
holomorphic case the full result (i.e., that the category of
$\V^G$-modules is equivalent to the category of modules over the
(twisted) Drinfeld double $D(G)$) is still not proved.

In this paper we suggest a new approach to the problem. The main idea
of this approach is {\bf not} using the structure theory of VOA's. In
the author's opinion, all the information about VOA's which is
relevant for this problem is encoded in the category of
representations of $\V$. For example, the pair $\V^G\subset \V$ can be
described in this way: as discussed in \cite{KO}, such a pair is the
same as an associative commutative algebra in the category $\C=\Rep
\V^G$ with some technical restrictions. Thus, if we know some basic
properties of $\C$ --- e.g., that this category is semisimple,
braided, and rigid --- then we can forget anything else about VOA's,
operator product expansions, etc. Instead, we use well-known tools for
working with braided tensor categories, such as graphical presentation
of morphisms.

Using this approach, in this paper we give an accurate proof of the
above conjecture; for simplicity, we only consider the case when all
the ``twists'', i.e. phase factors, are trivial. In this case, the
main result reads as follows.

\begin{theorem*}
  Let $\V$ be a VOA, $G$ -- a finite group of automorphisms of $\V$,
  and $\V^G$ -- the algebra of invariants. Assume that
\begin{enumerate}
\item $\V$ is ``holomorphic'', i.e. has a unique irreducible module,
  $\V$ itself, so that $\Rep \V=\vec$.
\item $\Rep \V^G$ is a semisimple  braided rigid balanced category
\item $\V$ has finite length as a $\V^G$-module
\item Certain cohomology class $\om\in H^3(G,\Ctimes)$ defined by $\V$
  is trivial
\end{enumerate}

Then the category $\Rep \V^G$ is equivalent to the category of modules
over $D(G,H)=\Cset[G]\ltimes\F(H)$ for some normal subgroup $H\subset
G$. If, in addition, we assume that $\Rep \V^G$ is modular, then $H=G$ so
$\Rep \V^G\simeq \Rep D(G)$.
\end{theorem*}

In this paper we assume that the reader is well familiar with braided
tensor categories and in particular, with the technique of using
graphs to prove  identities in such categories, developed by
Reshetikhin and Turaev. This can be found in many textbooks (see,
e.g., \cite{Ka}); we follow
the conventions of \cite{BK}. Conversely, knowledge of vertex operator
algebras and conformal field theory is not required: they do not even
appear in the paper except in this introduction. 

The paper also makes heavy use of results of \cite{KO}, so we suggest
that the reader keep a copy of that paper handy.

\section{Preliminaries}\label{s:preliminary}
  Throughout the paper, we denote by $\C$ a semisimple braided tensor
  category over $\Cset$, with simple objects $L_i, i\in I$ (``simple''
  always means ``non-zero simple''). As usual, we assume that the unit
  object is simple and denote the corresponding index in $I$ by $0$:
  $\one=L_0$. We assume that all spaces of morphisms are
  finite-dimensional and denote $\<V, W\>=\dim\Hom_\C(V, W)$; in
  particular, $\<L_i, V\>$ is the multiplicity of $L_i$ in $V$.
  
  As any abelian category, $\C$ is a module over the category $\vec$
  of finite-dimensional complex vector spaces, i.e. we have a natural
  functor of ``external tensor product'':
\begin{equation}
\tbox\colon \vec\times \C\to \C
\end{equation}
defined by $\Hom(L_i, V\tbox L)=V\otimes \Hom(L_i, L)$ (more formally:
$V\tbox L$ is the object representing the functor $F(M)=V\otimes
\Hom(M,L)$). This functor is bilinear and has natural associativity
properties:
\begin{equation}\label{e:assoc}
\begin{gathered}
  V\tbox (W\tbox L)=(V\ttt W)\tbox L\\
  V\tbox(L\ttt M)=(V\tbox L)\ttt M\\
  (V\tbox L)\ttt (W\tbox M)=(V\ttt W)\tbox (L\ttt M)
\end{gathered}
\end{equation}
(here $=$ means ``canonically isomorphic''). Abusing the language, we
will sometimes use $\ttt$ instead of $\tbox$.

Also, we denote by $G$ a compact group (e.g., $G$ can be finite) and
by $\Rep G$ the category of finite-dimensional complex representations of $G$.
This category is semisimple. We denote by $\Ghat$ the set of
isomorphism classes of simple $G$-modules, and for each $\la\in \Ghat$
we choose a representative $V_\la$.

We denote by $\C[G]$ the category whose  objects are pairs \\
$(M\in\C, \text{action of $G$ by automorphisms on }M)$. In particular,
each object of $\C$ can be considered as an object in $\C[G]$ by
letting $G$ act trivially. This category has the following properties,
proof of which is left as an exercise to the reader. 

\begin{enumerate}
\item $\C[G]$ is a semisimple rigid braided tensor category, with
  simple objects $V_\la\tbox L_i$.
\item Define the functor of $G$-invariants $\C[G]\to \C$ by
  $$
  \Hom_\C (L_i, M^G)= (\Hom_\C(L_i, M))^G=\Hom_{\C[G]}(L_i, M),
  $$
  i.e. $(V\tbox L_i)^G=V^G\tbox L_i$. Then this functor is exact,
  and one has canonical embedding
\begin{equation}\label{e:ginv}
X^G\ttt Y^G\injto (X\ttt Y)^G
\end{equation}

\item One can define canonical functor of ``exterior tensor product''
  $$
  \tbox\colon \Rep G\times \C[G]\to \C[G]
  $$
  which has the associativity properties \eqref{e:assoc}.
\end{enumerate}

\section{Untwisted sector}\label{s:main}
Throughout the paper,  we let $A$ be a $\C$-algebra (i.e. an object of $\C$
with a map $\mu \colon A\ttt A\to A$) as defined in \cite{KO}. We also
assume that $G$ acts on $A$ by multiplication-preserving automorphisms
and that $A^G=\one$ (recall that by axioms of $\C$-algebra, one has
canonical embedding $\one\injto A$ and the multiplicity $\<A,
\one\>=1$). In addition, we will also assume that $\C$ is rigid and
balanced, and that $A$ is rigid and satisfies $\th_A=\id$. 

Our main goal is to describe the category $\C$ in terms of the group
$G$ and the category $\Rep^0 A$ (see \cite{KO} for definitions). The
main motivation for this comes from the orbifold conformal field
theories, as explained in the introduction; in this situation,
$\C=\Rep \V^G$, and $A$ is $\V$ considered as a $\V^G$-module (cf.
\cite{KO}). We will freely use results and notation of \cite{KO}.

We start by describing the structure of $A$ as an object of
$\C$.

Define a functor $\Phi\colon \Rep G\to \C$ by
\begin{equation}\label{e:phi}
V\mapsto (V\tbox A)^G.
\end{equation}
In other words, if one writes
$$
A=\bigoplus_{\la\in \Ghat} V_\la \tbox M_\la
$$
for some $M_\la\in \C$, then $\Phi(V_\la^*)=M_\la$.

\begin{theorem}\label{t:phi(1)}
  $\Phi(\Cset)=\one$ and $\<\Phi(V_\la), \one\>=0$ for
  $V_\la\not\simeq\Cset$.
\end{theorem}
\begin{proof} Immediate from definitions. 
\end{proof}

Our next goal is to prove that under suitable conditions, $\Phi$
is a tensor functor. An impatient reader can find the final result as
\thref{t:main} below. We start by constructing a morphism $\Phi(V)\ttt
\Phi(W)\to \Phi (V\ttt W)$.
\begin{theorem}\label{t:J}
  Define the functorial morphism $J\colon \Phi(V)\ttt \Phi(W)\to \Phi
  (V\ttt W)$ as the following composition:
  $$
  (V\tbox A)^G\ttt (W\tbox A)^G\injto ((V\ttt W)\tbox(A\ttt
  A))^G\xxto{\mu} ((V\ttt W)\tbox A)^G
  $$
  \textup{(}we have used \eqref{e:assoc}, \eqref{e:ginv}\textup{)}.
  Then $J$ is compatible with associativity, commutativity,  unit, and
  balancing  morphisms in $\Rep G, \C$.
\end{theorem}
\begin{proof} Immediate from definitions. 
\end{proof}

This theorem allows one to define functorial morphisms $J\colon
\Phi(W_1)\ttt \dots \Phi(W_n)\to \Phi (W_1\ttt \dots\ttt W_n)$. For
future use, we explicitly write the functoriality property: for any
$f\colon W_i\to W'_i$, the following diagram is commutative
$$
\begin{CD}
  \Phi(W_1)\ttt \dots\ttt \Phi(W_n) @>{1\ttt\dots\ttt\Phi(f)
    \ttt\dots\ttt 1}>>
  \Phi(W_1)\ttt \dots\ttt \Phi(W'_i)\ttt\dots\ttt \Phi(W_n)\\
  @V{J}VV                          @V{J}VV\\
  \Phi(W_1\ttt \dots\ttt W_n) @>{\Phi(1\ttt\dots\ttt f \ttt\dots\ttt
    1)}>> \Phi(W_1\ttt \dots\ttt W'_i\ttt\dots\ttt W_n)
\end{CD}
$$

\begin{remark} We do not claim that $J$ is an isomorphism: in general,
  this is false.
\end{remark}

Now let us use rigidity of $A$.  We denote by $e_V\colon V^*\ttt V\to
\one, i_V\colon \one\to V\ttt V^*$ the canonical rigidity morphisms, and by
$\dim M$ dimension of an object $M\in \C$.

\begin{theorem} \label{t:v*}
  For every $V\in \Rep G$, the map
\begin{equation}
\tilde e\colon \Phi(V^*)\ttt \Phi(V)\xxto{J}\Phi(V^*\ttt
V)\xxto{\Phi(e)}\Phi(\Cset)=\one
\end{equation}
gives an isomorphism $\Phi(V^*)\simeq \Phi(V)^*$. In other words,
there exists a morphism $\tilde\imath\colon \one\to \Phi(V)\ttt
\Phi(V^*)$ such that $\tilde e, \tilde\imath$ satisfy the rigidity
axioms.
\end{theorem}
\begin{proof}
  It is easy to see that for any $X\in \C[G]$, the morphism
  $(X^*)^G\ttt X^G\injto (X^*\ttt X)^G\to \one^G=\one$ identifies
  $(X^*)^G\simeq (X^G)^*$ (here we have used \eqref{e:ginv} and
  functoriality of $X\mapsto X^G$). Combining this with the definition
  of $J$ and rigidity of $A$ we get the statement of the theorem.
\end{proof}

\begin{theorem}\label{t:Jinjective}
  For any $V, W\in \Rep G$, the morphism $J\colon \Phi(V)\ttt \Phi(W)\to
  \Phi(V\ttt W)$ is injective.
\end{theorem}

\begin{proof}
  Define the morphism $I\colon \Phi(V\ttt W)\to \Phi(V)\ttt \Phi(W)$
  by the graph shown in \firef{f:I}.
\begin{figure}
  \newfig{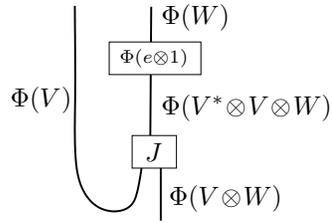}{%
           \rput(-1,1.3){$J$}
           \rput(-2.5,2){$\Phi(V)$}
           \rput[Bl](-0.8,0.6){$\Phi(V\sttt W)$}
           \rput[Bl](-0.9,1.8){$\Phi(V^*\sttt V\sttt W)$}
           \rput[Bl](-0.9,3){$\Phi(W)$}
           \rput(-1,2.55){$\scriptstyle\Phi(e\ttt 1)$}
}
\caption{Definition of $I$}\label{f:I}
\end{figure}

Then the composition $IJ\colon \Phi(V)\ttt \Phi(W)\to \Phi(V)\ttt
\Phi(W)$ can be rewritten as shown in \firef{f:IJ} and thus, $IJ=\id$
which proves that $J$ is injective.

\begin{figure}[h]
\begin{align*}
  &\newfig{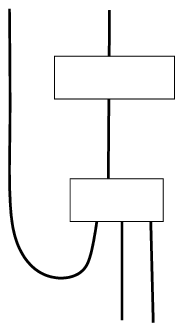}{
               \rput(-2.7,3){$\Phi(V)$}%
               \rput(-1.5,0.5){$\Phi(V)$}%
               \rput(-0.2,0.6){$\Phi(W)$}%
               \rput[Bl](-1, 2.2){$\Phi(V^*\sttt V\sttt W)$}%
               \rput[Bl](-1, 3.4){$\Phi(W)$}%
               \rput(-1.1,1.65){$J$}%
               \rput(-1.1,2.9){$\scriptstyle\Phi(e\ttt 1)$}%
           }\qquad\qquad\quad
   =\qquad\newfig{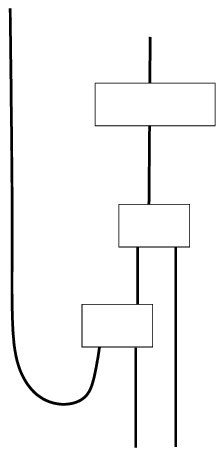}{%
               \rput(-2.9,3){$\Phi(V)$}%
               \rput(-1.7,0.5){$\Phi(V)$}%
               \rput(-0.2,0.6){$\Phi(W)$}%
               \rput[Bl](-0.95, 3.2){$\Phi(V^*\sttt V\sttt W)$}%
               \rput[Bl](-0.95, 4.4){$\Phi(W)$}%
               \rput(-1.3,1.65){$J$}%
               \rput(-1,2.65){$J$}%
               \rput(-0.9,3.9){$\scriptstyle\Phi(e\ttt 1)$}%
           }\\
  &\qquad =\qquad\newfig{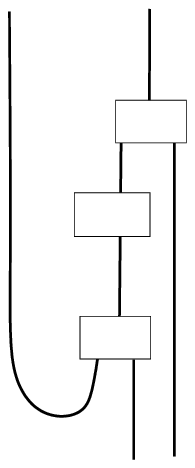}{%
               \rput(-2.7,3){$\Phi(V)$}%
               \rput(-1.35,0.5){$\Phi(V)$}%
               \rput(0,0.6){$\Phi(W)$}%
               \rput[Bl](-0.65, 4.4){$\Phi(W)$}%
               \rput(-1.1,1.65){$J$}%
               \rput(-1.2,3.3){$\one$}%
               \rput(-0.75,3.85){$J$}%
               \rput(-1.1,2.9){$\scriptstyle\Phi(e)$}%
          }\qquad
   =\qquad\newfig{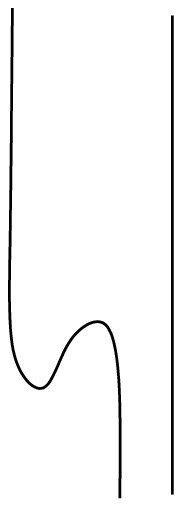}{
               \rput(-2.7,4){$\Phi(V)$}%
               \rput(0.05,4){$\Phi(W)$}%
           }
\end{align*}
\caption{Proof of $IJ=\id$}\label{f:IJ}
\end{figure}
\end{proof}

\begin{theorem}
  Let $V_\la$ be an irreducible representation of $G$.
\begin{enumerate}
\item $\Phi(V_\la)$ is either zero or a simple object in $\C$.
\item If $\la,\mu \in \Ghat, \la\ne \mu$ are such that $\Phi(V_\la)\ne
  0,\Phi(V_\mu)\ne 0$, then $\Phi(V_\la)\not \simeq \Phi(V_\mu)$.
\end{enumerate}
\end{theorem}
\begin{proof}
  Using \thref{t:phi(1)}, \thref{t:Jinjective}, we get
  $$
  \<\Phi(V_\la)\ttt \Phi(V_\mu^*), \one\>\le \<\Phi(V_\la\ttt
  V_\mu^*), \one\>=\de_{\la\mu}
  $$
  Thus, using $\Phi(V^*)\simeq\Phi(V)^*$ (see \thref{t:v*}) for
  $\la=\mu$ we get $\<\Phi(V_\la)\ttt \Phi(V_\la)^*, \one\>\leq 1$,
  which shows that $\Phi(V_\la)$ is either simple or zero, and for
  $\la\ne\mu$, we get $\<\Phi(V_\la)\ttt \Phi(V_\mu)^*, \one\>=0$.
\end{proof}

\begin{lemma}\label{l:lemma1}
  Let $V_\la$ be an irreducible representation of $G$
  such that $\Phi(V_\la)\ne 0$. Then the composition
  $$
  \one \xxto{\tilde\imath}\Phi(V_\la)\ttt \Phi(V_\la^*)
  \xxto{J}\Phi(V_\la\ttt V_\la^*)
  $$
  is equal to $\Phi(i_\la)$, where $i_\la$ is the canonical map
  $\Cset\to V_\la\ttt V_\la^*$.
\end{lemma}
\begin{remark}
  For $\Phi(V_\la)=0$, this lemma obviously fails.
\end{remark}

\begin{figure}
  \newfig{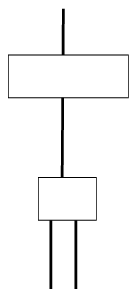}{
           \rput[Br](-1.35, 0.6){$\Phi(V_\la^*)$}
           \rput[Bl](-0.9, 0.6){$\Phi(V_\la\sttt V_\la^*)$}
           \rput[Bl](-1.05,1.8){$\Phi(V_\la^*\sttt V_\la\sttt V_\la^*)$}
           \rput(-1.1,2.55){$\scriptstyle\Phi(e\ttt 1)$}
           \rput(-1.1, 1.3){$J$}
           \rput[Bl](-1,3){$\Phi(V_\la^*)$}
            } 
\caption{Definition of $f\colon\Phi(V_\la^*)\ttt \Phi(V_\la\ttt
  V_\la^*)\to \Phi(V_\la^*)$}
\label{f:II}
\end{figure}

\begin{proof}
  Denote this composition by $\varphi$.  Since $\<\Phi(V_\la\ttt
  V_\la^*), \one\>=1$, one has $\varphi = c\Phi(i_\la)$ for some $c\in
  \Cset$. To find $c$, compose both $\varphi$ and $\Phi(i_\la)$ with
  the morphism $f\colon \Phi(V_\la^*)\ttt \Phi(V_\la\ttt V_\la^*)\to
  \Phi(V_\la^*)$ shown in \firef{f:II}. Arguing as in the proof of
  \thref{t:Jinjective}, we see that $f\varphi
  =f\Phi(i_\la)=\id_{\Phi(V_\la^*)}$.
\end{proof}

\begin{theorem}\label{t:Jisom}
  If  $\la\in \Ghat$ such that $\Phi(V_\la)\ne 0$, then
  $J\colon \Phi(V_\la)\ttt \Phi(W)\to \Phi(V_\la\ttt W)$ is an
  isomorphism for any $W\in \Rep G$.
\end{theorem}
\begin{proof}
  Let $I$ be as in the proof of \thref{t:Jinjective}. Let us calculate
  $JI$. Using \leref{l:lemma1} and functoriality of $J$, we can
  rewrite $JI$ as shown in \firef{f:JI}.

\begin{figure}
\begin{align*}
  \newfig{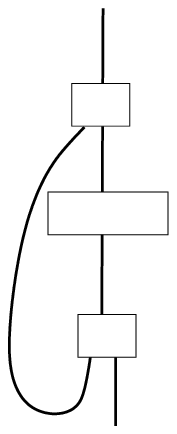}{           
           \rput(-1.1,1.3){$J$}
           \rput(-2.6,2){$\Phi(V)$}
           \rput[Bl](-0.9,0.6){$\Phi(V\sttt W)$}
           \rput[Bl](-1,1.8){$\Phi(V^*\sttt V\sttt W)$}
           \rput[Bl](-1,3){$\Phi(W)$}
           \rput(-1.1,2.55){$\scriptstyle\Phi(e\ttt 1)$}
           \rput(-1.1,3.65){$J$}
           \rput[Bl](-1,4.2){$\Phi(V\sttt W)$}
          }\qquad\qquad
  &=\newfig{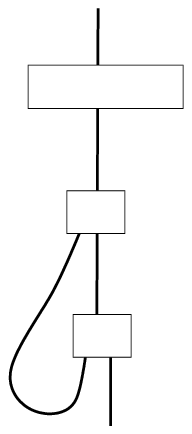}{
           \rput(-1.3,1.3){$J$}
           \rput(-2.4,2){$\Phi(V)$}
           \rput[Bl](-1.1,0.6){$\Phi(V\sttt W)$}
           \rput[Bl](-1.2,1.8){$\Phi(V^*\sttt V\sttt W)$}
           \rput[Bl](-1.2,3){$\Phi(V\sttt V^*\sttt V\sttt W)$}
           \rput(-1.3,3.8){$\scriptstyle\Phi(1\ttt e\ttt 1)$}
           \rput(-1.4,2.6){$J$}
           \rput[Bl](-1.2,4.3){$\Phi(V\sttt W)$}
          }\\[10pt]
  =\newfig{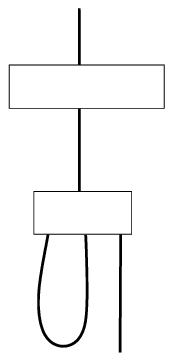}{
           \rput(-2.3,1){$\Phi(V)$}
           \rput[Bl](-0.9,0.6){$\Phi(V\sttt W)$}
           \rput(-1.3,1.8){$J$}
           \rput[Bl](-1.2,2.3){$\Phi(V\sttt V^*\sttt V\sttt W)$}
           \rput(-1.3,3.1){$\scriptstyle\Phi(1\ttt e\ttt 1)$}
           \rput[Bl](-1.2,3.6){$\Phi(V\sttt W)$}
           }\qquad\qquad\quad
 &=\newfig{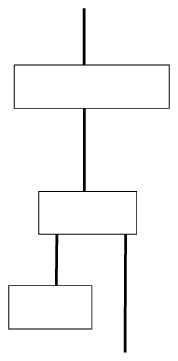}{
           \rput(-1.7,0.8){$\scriptstyle\Phi(i_\la)$}
           \rput[Br](-1.8,1.2){$\Phi(V\sttt V^*)$}
           \rput(-1.3,1.8){$J$}
           \rput[Bl](-0.9,0.6){$\Phi(V\sttt W)$}
           \rput[Bl](-1.2,2.3){$\Phi(V\sttt V^*\sttt V\sttt W)$}
           \rput(-1.3,3.1){$\scriptstyle\Phi(1\ttt e\ttt 1)$}
         }\qquad\qquad\quad
=\newfig{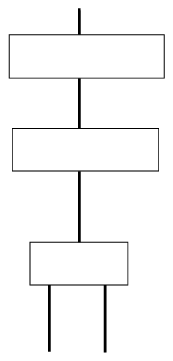}{
           \rput[Bl](-1,0.6){$\Phi(V\sttt W)$}
           \rput[Br](-1.8,0.6){$\one$}
           \rput(-1.3,1.3){$J$}
           \rput(-1.3,2.4){$\scriptstyle\Phi(i_V\ttt 1\ttt 1)$}
           \rput[Bl](-1.2,2.8){$\Phi(V\sttt V^*\sttt V\sttt W)$}
           \rput(-1.3,3.4){$\scriptstyle\Phi(1\ttt e\ttt 1)$}
          }        
\end{align*}
\caption{Proof of $JI=\id$}
\label{f:JI}
\end{figure}

Thus, $JI=\id$. On the other hand, it was proved in
\thref{t:Jinjective} that $IJ=\id$.
\end{proof}

Now, let us assume that the action of $G$ on $A$ is faithful, that is,
every $g\in G, g\ne 1$ acts on $A$ by a non-trivial automorphism.

\begin{theorem}
  If the action of $G$ is faithful, then
  $\Phi(V_\la)\ne 0$ for any $\la\in \Ghat$.
\end{theorem}

\begin{proof}
  Let $I=\{\la\in \Ghat\st \Phi(V_\la)\ne 0\}$. Then, by \thref{t:v*},
  $I$ is closed under duality (as usual, we denote by $\la^*$ the
  class of representation $V_\la^*$). Denote by $N_{\la\mu}^\nu$ the
  multiplicities in the tensor product decomposition: $V_\la\ttt
  V_\mu\simeq \sum N_{\la\mu}^\nu V_\nu$. Then
\begin{equation}\label{e:Nne0}
N_{\la\mu}^\nu=0\qquad \text{ if }\la,\nu\in I, \mu\notin I.
\end{equation}
Indeed: if $N_{\la\mu}^\nu\ne 0$, then there is an embedding
$V_\nu\subset V_\la\ttt V_\mu$, which gives
$$
\Phi(V_\nu)\subset \Phi(V_\la\ttt V_\mu)\simeq \Phi(V_\la)\ttt
\Phi(V_\mu) \simeq 0
$$
(by \thref{t:Jisom}), which contradicts $\Phi(V_\nu)\ne 0$.

By using $N_{\la\mu}^\nu=N_{\la\nu^*}^{\mu^*}$, we can rewrite
\eqref{e:Nne0} as
\begin{equation}\label{e:Nne0-2}
N_{\la\mu}^\nu=0\qquad \text{if }\la,\mu\in I, \nu\notin I.
\end{equation}
Let $\A$ be the full subcategory in $\Rep G$ generated (as an abelian
category) by $V_\la, \la\in I$. Then it follows from \eqref{e:Nne0-2}
that $\A$ is closed under tensor product; it is also closed under
duality (by \thref{t:v*}) and contains $\Cset$ (by \thref{t:phi(1)}).
By the usual reconstruction theorems, this means that $\A$ is the
category of representations of some group $H$ which is a quotient of
$G$. But by assumption, the action of $G$ on $A$ is faithful, which
means that the action of $G$ on $\bigoplus_{\la\in I}V_\la$ is
faithful. Thus, $H=G$, $I=\Ghat$.
\end{proof}

Combining all the results above, we get the main theorem of this
section.
\begin{theorem}\label{t:main}
  Let $A$ be a rigid $\C$-algebra, $\th_A=1$, and $G$ --- a compact
  group acting faithfully on $A$. Then
\begin{enumerate}
\item
\begin{equation}\label{e:decomposition}
A\simeq \bigoplus_{\la\in \Ghat}V_\la \tbox M_\la
\end{equation}
where $M_\la=\Phi(V_\la^*)$ are non-zero, simple, and
$M_\la\not\simeq M_\mu$ for $\la\ne \mu$.
\item Let $\C_1$ be the full subcategory in $\C$ generated as an
  abelian category by $M_\la, \la\in\Ghat$. Then
  $\C_1$ is a symmetric tensor subcategory in $\C$, and the functor
  $\Phi\colon \Rep G\to \C_1$ defined by \eqref{e:phi} is an
  equivalence of tensor categories.
\end{enumerate}
\end{theorem}

\begin{corollary}\label{c:gfinite}
  $G$ is finite.
\end{corollary}
\begin{proof}
Immediate from \eqref{e:decomposition} and finite-dimensionality of
spaces of morphisms in $\C$. 
\end{proof}

\begin{corollary}\label{c:R^2}
  $\check R_{AA}^2=\id$.
\end{corollary}

\begin{example}\label{e:repG}
  Let $G$ be a finite group, $\C=\Rep G$, $A=\F(G)$ --- the algebra of
  functions on $G$, with the usual (pointwise) multiplication. Formula
  $g\de_h=\de_{gh}$ makes $A$ a $G$-module and thus, an object of
  $\C$. It is trivial to show (see \cite{KO}) that $A$ is a rigid
  $\C$-algebra, and $\Rep A=\vec $.
  
  We also have another action of $G$ on $A$, by
  $\pi_g\de_h=\de_{hg^{-1}}$. This commutes with previously defined
  and thus, defines an action of $G$ by automorphisms on $A$
  considered as a $\C$-algebra. In this case, the functor $\Phi$ is an
  equivalence of categories $\Rep G\simeq \C$, and the decomposition
  \eqref{e:decomposition} becomes the standard decomposition
  $$
  \F(G)=\bigoplus_{\la\in \Ghat} V_\la\ttt V_\la^*.
  $$

\end{example}

Let us return to the general case.
\begin{theorem}\label{t:a=f(G)}
  Under the assumptions of \thref{t:main}, consider $A\in
  \C_1\subset \C$ as an object of $\Rep G$ using equivalence $\Phi$.
  Then $A\simeq\F(G)$ with multiplication, structure of $G$-module and
  action of $G$ by automorphisms as defined in \exref{e:repG}.
\end{theorem}
\begin{proof}
  It is immediate from \thref{t:main} that $A$ lies in $\C_1\simeq
  \Rep(G)$ and that as an object of $\Rep (G)$, $A=\bigoplus_{\la\in
    \Ghat} V_\la \tbox V_\la^*$. The structure of $G$-module is
  determined by the action on the second factor, and the action of $G$
  by automorphism is defined by the action on the first factor.
  
  On the other hand, it is shown in  \cite{KO} that any
  algebra in category $\Rep G$ must be of the form $A=\F(G/H)$ for
  some subgroup $H$. Combining these statements, we see that $H=\{1\},
  A=\F(G)$.
\end{proof}

\begin{corollary}\label{c:autom}
  $G$ is the group of all automorphisms of $A$ as a $\C$-algebra.
\end{corollary}

\begin{corollary}\label{c:dima}
  $\dim A=|G|$.
\end{corollary}

Of course, in general case $\C$ can be (and usually is) larger than
$\Rep G$. A very important example is when $\C$ is the category of
representation of $D(G)$ (the Drinfeld double of $G$), and
$A=\F(G)\in\Rep G\subset \Rep D(G)$ with the same action of $G$ as in
\exref{e:repG}. This example plays an important role in what follows;
it is discussed in detail in \seref{s:D(G)}.

\begin{example}\label{x:simplecurr}
Let $G$ be commutative. Then all its irreducible representations are
one-dimensional, and it immediately follows from \thref{t:main} that
$$
A\simeq \bigoplus_{\la\in \Ghat} M_\la
$$
where $M_\la$ are non-zero, simple, pairwise non-isomorphic objects
in $\C$ and $M_\la\ttt M_\mu\simeq M_{\la\mu}$. This case is well
studied in numerous papers under the name ``simple currents
extensions''; a review of known results can be found, e.g., in
\cite{CFT} and \cite{FS}.
\end{example}

\section{Example: $D(G)$}\label{s:D(G)} 
Let $D(G)=\Cset[G]\ltimes \F(G)$ be the Drinfeld's double of the
finite group $G$ (see, e.g., \cite{BK}), and let $\C$ be the category
of finite-dimensional complex $D(G)$-modules. As is well known, this
category is equivalent to the category of $G$-equivariant vector
bundles on $G$. An object of this category is a complex vector space
$V$ with an action of $G$ and with a $G$-grading: $V=\bigoplus_{g\in
  G} V_g$ such that $gV_x\subset V_{gxg^{-1}}$, or, equivalently,
$\wt(gv)=g \wt(v)g^{-1}$, where $\wt(v)\in G$ denotes weight of a
homogeneous $v\in V$. The tensor product in $\C$ is the usual tensor
product, with $\wt(v\ttt w)=\wt(v)\wt(w)$.  The braiding in $\C$ is
given by $\check R=PR$, where
\begin{equation}\label{e:rdg}
R=\sum_{g\in G}\de_g\ttt g\in D(G)\ttt D(G).
\end{equation}
In other words, if $v,w$ are homogeneous vectors in $V, W$
respectively, then
$$
\check R(v\ttt w)=gw\ttt v, \quad g=wt(v).
$$

It is also known that $\Rep D(G)$ is semisimple, and the set of
isomorphism classes of simple objects is $(g,\pi)/G$ where $g\in G,
\pi$ -- an irreducible representation of the centralizer $Z(g)=\{h\in
G\st hg=gh\}$, and $G$ acts on the set of pairs $(g,\pi)$ by
$h(g,\pi)=(hgh^{-1}, \pi\circ h^{-1})$. We will denote the
corresponding representation of $D(G)$ by $V_{g,\pi}$.

Let $A=\F(G)$; consider it as object of $\C$ by endowing it with the
standard action of $G$ (same as in \exref{e:repG}), and by letting
$wt(v)=1$ for all $v\in A$.
\begin{lemma}
  $A$ is a rigid $\C$-algebra, with $\th_A=1$. 
\end{lemma}
The proof is straightforward.

As in \exref{e:repG}, we also have an action of $G$ by automorphisms
on $A$ defined by $\pi_g\de_h=\de_{hg^{-1}}$, and
$A^G=\one=\Cset\bigl(\sum_{h\in G}\de_h\bigr)$.
 
We remind the reader that for every $\C$-algebra $A$, one can define
the category $\Rep A$ of $A$-modules and two natural functors $F\colon
\C\to \Rep A, G\colon \Rep A\to \C$ (see \cite{KO}). The category
$\Rep A$ is a tensor category; we denote tensor product in $\Rep A$ by
$\tta$.

\begin{theorem}\label{t:ggraded}
  The category $\Rep A$ is equivalent to the category $G\vec$ of
  $G$-graded vector spaces. Under this equivalence, the functor $\tta$
  becomes the usual tensor product of vector spaces with the grading
  given by $\wt(v\ttt w)=\wt(v)\wt(w)$, and the functors $F, G$ are
  given by
\begin{alignat*}{2}
  F(V)&=V & &\text{\rm forgetting the action of $G$ but keeping the grading}\\
  G(V)&=\F(G)\ttt V &\quad&\text{\rm with grading given by
    $\wt(\de_g\ttt v)=g\wt(v)g^{-1}$}
\end{alignat*}
\end{theorem}
\begin{proof}
  It is immediate from the definition that $\Rep A$ is the category of
  $G$-modules with $G\times G$-grading such that
  $$
  \wt(gv)=(gxg^{-1},gy)\quad\text{if }\wt(v)=(x,y)\in G\times G
  $$
  Indeed, the action of $G$ and the first component of the grading
  define $V$ as an object of $\C$, and the second component of the
  grading defines the action of $A$: if $v\in V$ is homogeneous, then
  $$
  \de_gv=\begin{cases} v,& \wt(v)=(\cdot, g)\\
    0, &\text{otherwise}
        \end{cases}
        $$
        
        Define a new $G\times G$-grading $\twt$ by
        $$
        \twt(v)=(y^{-1}xy,y)\quad\text{if }\wt(v)=(x,y)
        $$
        (which implies that $\wt(v)=(y\tilde x y^{-1},y)$ if
        $\twt(v)=(\tilde x, y)$). Then
        $$
        \twt(gv)=(\tilde x, gy)\quad\text{if }\twt(v)=(\tilde x,
        y).
        $$
        From this it immediately follows that the functor $\Rep
        A\to G\vec$ given by
\begin{equation}
V\mapsto \{v\in V\st \wt(v)=(\cdot,1)\}=\{v\in V\st \twt(v)=(\cdot,1)\}
\end{equation}
considered with $G$-grading given by the first component of $\wt$, is
an equivalence of categories. The remaining statements of the theorem
are straightforward.
\end{proof}

\begin{corollary}
  The set of simple objects in $\Rep A$ is $G$.
\end{corollary}
For future use we give description of the corresponding simple objects
$X_g$ in terms of $G$-graded vector spaces and in terms of $G\times
G$-graded $G$-modules. As a $G$-graded vector space,
$$
(X_g)_h=\begin{cases} \Cset, &g=h\\
  0,& \text{otherwise}
        \end{cases}
$$
As a bi-graded $G$-module, $X_g$ is given by
\begin{equation}\label{e:X_g}
  X_g=\bigoplus_{y^{-1}xy=g}\Cset e_{x,y}
\end{equation}
with $\wt(e_{x,y})=(x,y)$, the action of $G$ and $A$ given
by
\begin{align*}
  h e_{x,y}&=e_{hxh^{-1}, hy}\\
  \mu(\de_h\ttt e_{x,y})&=\begin{cases}
    e_{x,y}, &h=y\\
    0, &\text{otherwise}
                       \end{cases}
\end{align*}
Note that the first component of $\wt(v), v\in X_g$ is supported on
the conjugacy class of $g$.

This description immediately implies the following result: 
\begin{equation}\label{e:g(xg)}
G(X_g)=V_{g,\F(Z(g))}=
   \bigoplus_{\pi\in \widehat{Z(g)}} \pi\ttt V_{g,\pi}
\end{equation}
where $\pi$ is the representation space, considered with trivial
action of $D(G)$.

\begin{theorem}
  $\Rep^0 A=\vec$ which is considered as a subcategory in $G\vec$
  consisting of spaces with grading identically equal to 1.
\end{theorem}
\begin{proof}
  Let $V\in \Rep A$; for now, we consider $V$ as a $G\times G$-graded
  $G$-module, as in the proof of \thref{t:ggraded}. Then explicit
  calculation shows that
  $$
  \check R_{VA}\check R_{AV}(\de_g\ttt v)=\de_{xg}\ttt v
  \quad\text{if }\wt(v)=(x, \cdot).
  $$
  Therefore,
  $$
  \mu_V\check R_{VA}\check R_{AV}(\de_g\ttt v)=
\begin{cases}
  v,& xg=y\\
  0,&\text{otherwise}
\end{cases}
$$
where $\wt(v)=(x,y)$.  Comparing it with the usual formula for
action of $A$,
$$
\mu_V(\de_g\ttt v)=
\begin{cases}
  v,& g=y\\
  0,&\text{otherwise}
\end{cases}
$$
we see that $\mu_V=\mu_V\check R_{VA}\check R_{AV}$ iff
$\wt(v)=(1,\cdot)$ for all $v\in V$.
\end{proof}

For future use, we give here two more results about $D(G)$. First,
define the map 
\begin{equation}\label{e:tau}
\begin{aligned}
\tau\colon D(G)&\to D(G)\\
           g\de_h&\mapsto g\de_{h^{-1}}
\end{aligned}
\end{equation}         
Then it is trivial to check the following properties. 

\begin{lemma}\label{l:tau}
\begin{enumerate}
\item
 $\tau$ is an algebra automorphism. 
\item $\tau$ is coalgebra anti-automorphism: $\De \tau (a)=
  \tau\ttt\tau (\De^{op} a)$, $\tau\circ S=S\circ \tau$, where $S$ is
  the antipode in $D(G)$. 
\item $\bigl(\tau\ttt \tau\bigr) (R)=R^{-1}$.
\end{enumerate}
\end{lemma}
Thus, if we define for a representation $V$ of $D(G)$ a new
representation $V^{\tau}$ which coincides with $V$ as a vector space
but with the action of $D(G)$ twisted by $\tau$:
$\pi_{V^\tau}(a)=\pi_V(\tau a)$, then one has canonical isomorphisms
$(V\ttt W)^\tau=W^\tau\ttt V^\tau, (V^*)^\tau=(V^\tau)^*$ and thus,
$\tau$ gives an equivalence
\begin{equation}\label{e:repop}
\tau\colon (\Rep D(G))^{op}\isoto \Rep D(G),
\end{equation}
where $(\Rep D(G))^{op}$ coincides with $\Rep D(G)$ as an abelian
category but has tensor product and braiding defined by $V\ttt^{op}
W=W\ttt V, R^{op}=R^{-1}$.

Second, note that  $D(G)$ can be generalized as follows.
Let $H\subset G$ be a normal subgroup.
Define
\begin{equation}\label{e:dhg}
D(G,H)=\Cset[G]\ltimes \F(H).
\end{equation}
One easily sees that it is a quotient of $D(G)$: $D(G,H)=D(G)/J_H$,
where $J_H$ is the ideal generated by $\de_g, g\in G\setminus H$ (this
is a Hopf ideal).  For $H=\{1\}$, we get $D(G,H)=\Cset[G]$; for $H=G$,
$D(G,H)=D(G)$. One also easily sees that $\Rep D(G,H)$ is the
subcategory of $\Rep D(G)$ consisting of representations such that
$\wt(v)\in H$.

\begin{theorem}\label{t:rephg}
  Let $\A\subset \Rep D(G)$ be a full subcategory containing $\Rep G$
  \textup{(}which we consider as a subcategory in $\Rep D(G)$ as
  before\textup{)}, and
  closed under duality, tensor product, and taking  sub-objects. Then
  $\A=\Rep D(G,H)$ for some normal subgroup $H\subset G$.
\end{theorem}
\begin{proof}
  The proof is based on the following easily proved lemma.

\begin{lemma}
  Let $g_1,g_2\in G$ and let $\pi_1, \pi'$ be irreducible
  representations of $Z(g_1)$, $Z(g_1g_2)$ respectively. Then there
  exists $\pi_2$ --- an irreducible representation of $Z(g_2)$ such
  that $\<V_{g_1,\pi_1}\ttt V_{g_2, \pi_2}, V_{g_1g_2, \pi'}\>\ne 0$.
\end{lemma}

Using this lemma with $g_2=1$, we see that if $V_{g_1,\pi_1}\in \A$,
for some $\pi_1\in \widehat{Z(g_1)}$, then for all $\pi\in
\widehat{Z(g_1)}$, $V_{g_1,\pi}\in \A$. Using this lemma again, we see
that the set $H=\{g\in G\st V_{g,\pi}\in \A\}$ is closed under product.
Since $V_{g,\pi}=V_{hgh^{-1}, \pi\circ h^{-1}}$, this set is also
invariant under conjugation. Thus, $H$ is a normal subgroup in $G$.

\end{proof}

We  also need the following theorem.

\begin{theorem}\label{t:modularh=g}
$\Rep D(G,H)$ is modular iff $H=G$.
\end{theorem}
\begin{proof}
  It is well known that  $\Rep D(G)$ is modular; thus, let us prove that
  if $H\ne G$ then $\Rep D(G,H)$ is not modular. As discussed above,
  simple objects in $\Rep D(H,G)$ are $V_{g,\pi}, g\in H$. Let $\pi$
  be a formal linear combination of irreducible representations of $G$
  such that $\tr_\pi(h)=0$ for all $h\in H$; such a $\pi$ always
  exists if $H\ne G$. Then it follows from explicit formulas for $s$
  (see, e.g., \cite{BK}) that 
$$
s_{(1,\pi), (h,\pi')}=0
$$
for all $h\in H, \pi'\in \widehat{Z(h)}$. Thus, $s$ is singular.
\end{proof}

\section{Twisted modules}\label{s:twisted1}
As before, we let $A$ be a rigid $\C$-algebra with $\th_A=\id$, and
$G$ --- a compact group acting faithfully on $G$ by automorphisms (by
\coref{c:gfinite}, this implies that $G$ is finite). For every $g\in
G$, we will denote the corresponding automorphism of $A$ by $\pi_g$.
We use the same notation $\Rep A, \Rep^0 A, \mu_V\colon A\ttt V\to V$
and functors $F\colon \C\to \Rep A, G\colon \Rep A\to \C$ as in
\cite{KO}. We will also use the same conventions in figures as in
\cite{KO}, representing $A$ by dashed line.

From now on, we will also assume that $A$ is such that $\Rep^0 A$ has
a unique simple module, $A$ itself; thus, $\Rep^0 A\simeq \vec$. This
corresponds to ``holomorphic'' case in conformal field theory; for
this reason, we will call such $A$ ``holomorphic''.

\begin{figure}
  \newfig{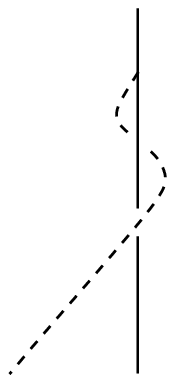}{\rput(-0.9,3.8){$X$}}
  =\newfig{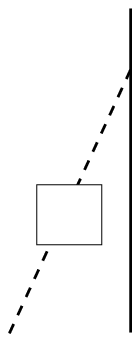}{\rput(-1,1.6){$\pi_g^{-1}$}
                          \rput(-0.6,3.4){$X$}}
\caption{Definition of $g$-twisted module} \label{f:twisted}
\end{figure}

\begin{definition}\label{d:twisted}
  Let $g\in G$. A module $X\in \Rep A$ is called $g$-twisted if
  $$
  \mu\check R^2=\mu\circ (\pi_{g^{-1}}\ttt \id)\colon A\ttt V\to V
  $$
  (see \firef{f:twisted}).
\end{definition}
In particular, $X$ is $1$-twisted iff $X\in \Rep^0A$. This definition
is, of course, nothing but rewriting in our language of the definition
given in \cite{DVVV}.

\begin{example}
  In the situation of \seref{s:D(G)}, the simple module $X_g$ is
  $g$-twisted. Indeed,
  $$
  \mu\check R^2(\de_h \ttt e_{x,y})=
           \begin{cases}
             e_{x,y}, & xh=y\\
             0, &\text{otherwise}
           \end{cases}
  $$
  But for $X_g$, $y^{-1}xy=g$, which trivially implies that $xh=y\iff
  hg=y$.  Thus, $\mu\check R^2(\de_h \ttt e_{x,y})=\mu(\de_{hg} \ttt
  e_{x,y})$.
\end{example}

The key result of this section is the following theorem.

\begin{theorem}\label{t:twisted1}
  Assume that $\Rep^0 A=\vec$. Then every simple object $X_i\in \Rep
  A$ is $g$-twisted for some $g=g(X_i)\in G$.
\end{theorem}
\begin{proof}
  The proof is based on the following lemma.
\begin{lemma}\label{l:projector}
  If $X\in \Rep A$ is simple and $A$ is holomorphic, then
  $$
  \frac{1}{(\dim A)^2}\fig{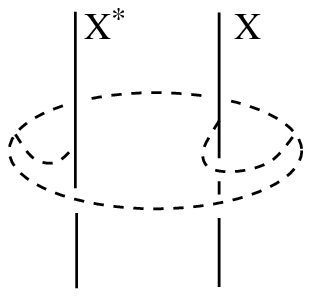}=\frac{1}{\dim
    X}\fig{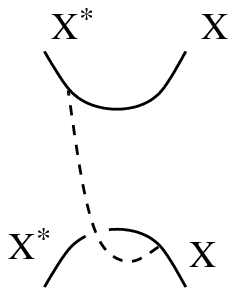}
  $$
\end{lemma}

Note that $\dim X_i$ is non-zero in any semisimple rigid category in
which $X^{**}\simeq X$ (see, e.g., \cite[Section 2.4]{BK}). In
particular, this implies $\dim_{\Rep A}X\ne 0$ and thus, $\dim
X=(\dim A)(\dim_{\Rep A}X)\ne 0$ (see \cite[Theorem 3.5]{KO}). 

\begin{proof}
  Let us rewrite the left hand side as shown in \firef{f:proj2}. Using
  Lemma~1.15, Lemma~5.3 from \cite{KO}, we see that the left hand side
  is the composition
  $$
  X^*\ttt X\to X^*\tta X\xxto{P} (X^*\tta X)^0
  $$
  where $P$ is the projector $\Rep A\to \Rep^0 A$, and for $V\in
  \Rep A$, $V^0=P(V)$ is the maximal sub-object of $V$ which lies in
  $\Rep^0A$. But if $A$ is holomorphic, then the only simple object in
  $\Rep^0A$ is $A$ itself, which is the unit object in $\Rep A$.  It
  appears in $X^*\tta X$ with multiplicity one, and the right-hand
  side is exactly the projection of $X^*\tta X$ on $A\subset X^*\tta
  X$.
\begin{figure}[h]
  $$
  \frac{1}{(\dim A)^2}\fig{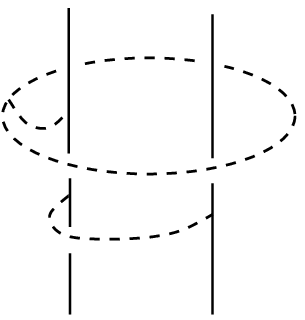}
  $$
\caption{}\label{f:proj2}
\end{figure}

\end{proof}

For every $X\in \Rep A$, define morphism $T_X\colon A\to A$ by the
following graph:
\begin{equation}\label{e:phiv}
T_X=\fig{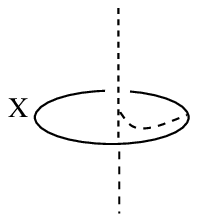}
\end{equation}

\begin{lemma}\label{l:phiautom}
  If $X\in \Rep A$ is simple, then $\frac{1}{\dim X}T_X$ is an
  automorphism of $A$ as a $\C$-algebra.
\end{lemma}

\begin{proof}
  Let us calculate $\mu\circ(T_X\ttt T_X)$. Using \leref{l:projector}
  we can rewrite the graph defining $\mu\circ(T_X\ttt T_X)$ as shown
  in \firef{f:TT}. (Note that we need to use $\check R^2_{AA}=1$
  (\coref{c:R^2}) to move the ``ring'' through $A$; the last step also
  uses Lemma~1.10 from \cite{KO}.)  This shows that $\mu\circ(T_X\ttt
  T_X)=(\dim X)T_X\circ \mu$. In other words, $\frac{1}{\dim X}T_X$ is
  a morphism of $\C$-algebras. But it easily follows from
  \thref{t:a=f(G)} that every such morphism is either zero or
  invertible. Restricting $T_X$ to $\one\subset A$, we see that $T_X$
  is non-zero; thus, $\frac{1}{\dim X}T_X$ is an automorphism.
\end{proof}

\begin{figure}
\begin{align*}
  &\fig{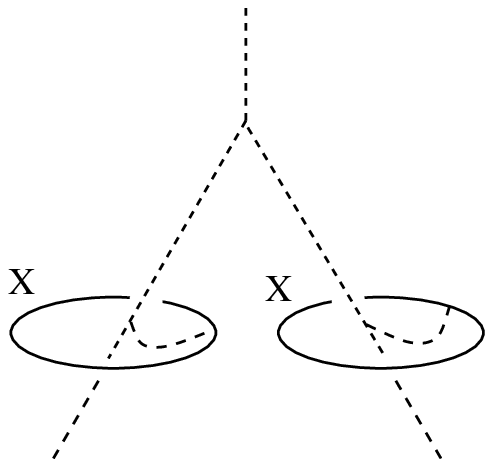}=\fig{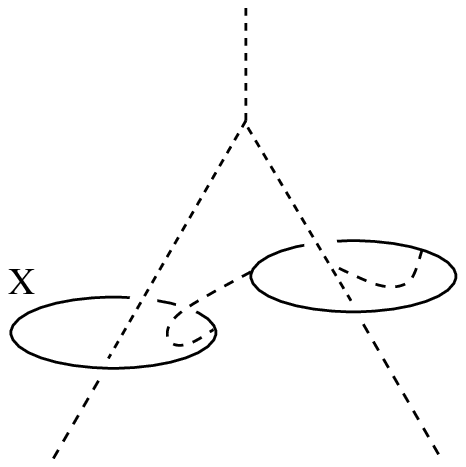}\\
  &=\frac{\dim X}{(\dim A)^2}\fig{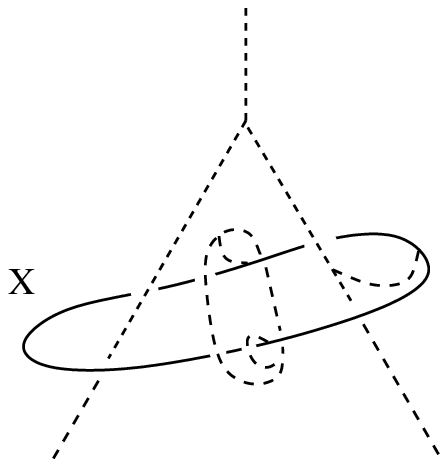}
  =\frac{\dim X}{(\dim A)^2}\fig{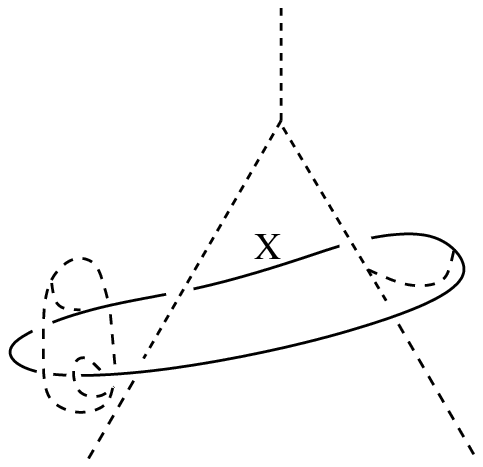}\\
  &=\frac{\dim X}{(\dim A)^2}\fig{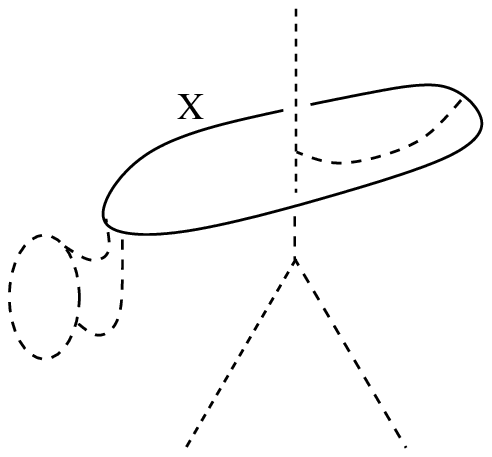} =(\dim X) T_X \circ \mu
\end{align*}
\caption{Proof of \leref{l:phiautom}} 
\label{f:TT}
\end{figure}

\begin{lemma}\label{l:phitwist}
  $$
  \frac{1}{\dim X}
    \newfig{twistedb.eps}{\rput(-1,1.6){$T_X$}
                          \rput(-0.6,3.4){$X$}}
=\fig{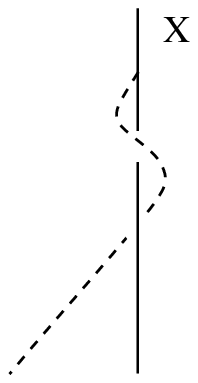}
  $$
\end{lemma}
\begin{proof}
  Using \leref{l:projector} and $\check R_{AA}^2=\id$, we can rewrite
  the left hand side as shown in \firef{f:phitwist}.

\begin{figure}
\begin{align*}
  &\frac{1}{\dim X}\fig{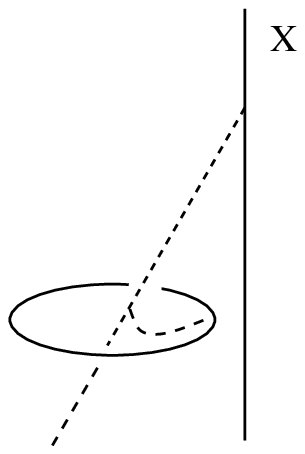}
   =\frac{1}{\dim X}\fig{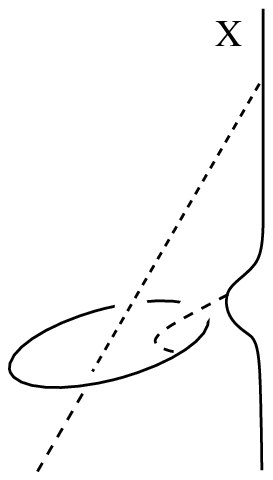}
   =\frac{1}{(\dim A)^2}\fig{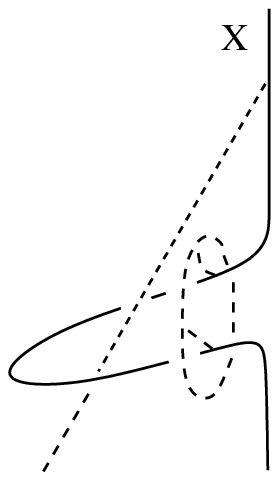}\\
  &=\frac{1}{(\dim A)^2}\fig{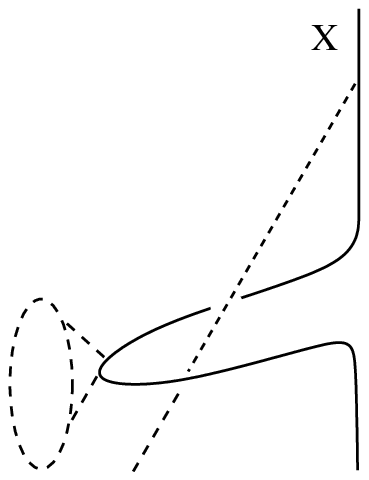}=\fig{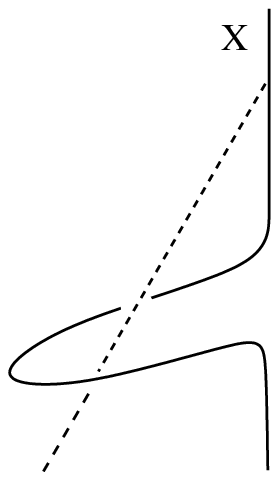}
\end{align*}
\caption{Proof of \leref{l:phitwist}}
\label{f:phitwist}
\end{figure}
\end{proof}

The statement of the theorem easily follows from these two lemmas.
Indeed, combining \leref{l:phiautom} with \coref{c:autom} we see that
$\frac{1}{\dim X}T_X=\pi_g$ for some $g\in G$. On the other hand,
\leref{l:phitwist} gives $\mu\circ \check
R^{-2}=\mu\circ\bigl(\frac{1}{\dim X}T_X\ttt
\id\bigr)=\mu\circ(\pi_g\ttt\id)$ and thus, $\mu\circ \check
R^{2}=\mu\circ(\pi^{-1}_g\ttt\id)$. This completes the proof of
\thref{t:twisted1}
\end{proof}
Let us study some properties of the correspondence $X\mapsto g(X)$.

\begin{theorem} \label{t:x_g}
  Let $X\in\Rep A$ be simple. Then
\begin{enumerate}
\item For $h\in G$, let $X^h\in \Rep A$ coincide with $X$ as an object
  of $\C$, but with a twisted action of $A:
  \mu_{X^h}=\mu_X\circ(\pi_h\ttt\id)$. Then $g(X^h)=h^{-1}g(X)h$.
\item $g(X)=1$ iff $X\simeq A$.
\item $g(X^*)=g(X)^{-1}$.
\item If $X, Y\in \Rep A$ are simple then so is $X\tta Y$ and $g(X\tta
  Y)=g(X) g(Y)$.
\item If $X, Y\in \Rep A$ are simple and non-isomorphic, then $g(X)\ne
  g(Y)$.
\end{enumerate}
\end{theorem}
\begin{proof}
  (1) is immediate from the definitions. (2) follows from the fact
  that the only simple object in $\Rep^0A$ is $A$. To prove (3) and
  (4), we will use the following lemma which easily follows from
  Lemma~1.15 in \cite{KO}.

\begin{lemma}\label{l:txy}
  $T_{X\tta Y}=\frac{1}{\dim A} T_XT_Y$.
\end{lemma}

This lemma implies that $T_{X\tta X^*}=c \pi_{g(X)g(X^*)}$ for some
$c\in \Cset$. On the other hand, $X\tta X^*\simeq A\oplus \sum_{i\ne
  0}N_iX_i$, and $T_{X\tta X^*}=\pi_1 +\sum c_i\pi_{g(X_i)}$. Since, by
(2), $g(X_i)\ne 1$ for $i\ne 0$ and the operators $\pi_g$ are linearly
independent, we see that these two expressions can be equal only if
$X\tta X^*\simeq A, g(X)g(X^*)=1$.

To prove (4), note that $X^*\tta (X\tta Y)\simeq (X^*\tta X)\tta
Y\simeq A\tta Y\simeq Y$ is simple, which immediately implies that
$X\tta Y$ is simple. The identity $g(X\tta Y)=g(X)g(Y)$ immediately
follows from \leref{l:txy}.

Finally, to prove (5) note that (3) and (4) imply $g(X)=g(Y)\iff
g(X^*)g(Y)=g(X^*\tta Y)=1$. By (2), this is only possible if $X^*\tta
Y\simeq A$.
\end{proof} 

\begin{corollary}\label{c:H}
  The map $X\mapsto g(X)$ is a bijection between the set of
  isomorphism classes of simple objects in $\Rep A$ and some normal
  subgroup $H\subset G$.
\end{corollary}

We will denote the unique simple $g$-twisted object $X\in \Rep A$ by
$X_g$; then \thref{t:x_g} implies that
\begin{equation}\label{e:xgh}
X_{g}\tta X_h\simeq X_{gh}.
\end{equation}
Combining this with multiplicativity of dimension and
\cite[Theorem 3.5]{KO}, we see that
\begin{equation*}
g\mapsto \dim_A X=\frac{\dim X}{\dim A}
\end{equation*}
is a character of the group $H\subset G$. Since $H$ is a finite
group, this immediately implies the following result.

\begin{lemma}
  For any simple $X\in \Rep A$, $\left|\frac{\dim X}{\dim
      A}\right|=1$.
\end{lemma}

In particular, if $\dim X\ge 0$ (which happens if $\C$ is a unitary
category in the sense of Turaev), then this implies $\dim X_g=\dim
A=|G|$.

\section{Twisted sectors}\label{s:twisted2}
Now that we have a description of irreducible objects in $\Rep A$ in
terms of the group $G$, we can move on to our ultimate goal:
description of irreducible objects in $\C$. As before, we assume that
$A$ is rigid, $\th_A=\id$ and holomorphic ($\Rep^0 A=\vec$), and $G$
is a finite group acting faithfully on $A$.

It is immediate from the identity $\<M, G(X)\>_\C=\<F(M), X\>_{\Rep
  A}$ (see \cite{KO}) that every simple $L_i\in \C$ appears in the
decomposition of some $X_g\in \Rep A$ (considered as an object of
$\C$).  Thus, our first goal is to study the decomposition
\begin{equation}\label{e:x_g}
X_g\simeq \bigoplus N_{g,i}L_i.
\end{equation}
Note that it immediately follows from \thref{t:x_g} that as an object
of $\C$, $X_g\simeq X_{hgh^{-1}}$; thus, the multiplicities $N_{g,i}$
only depend on the conjugacy class of $g$.

Our strategy in studying  decomposition \eqref{e:x_g} is parallel
to the approach taken in \seref{s:main} to study the decomposition of
$A=X_1$. However, instead of the functor $\Phi\colon \Rep G\to \C$
which was defined using $A$, we will define functor $\Phi\colon \Rep
D(G)\to \C$ using $\tA=\bigoplus_{g\in G}X_g$, where $D(G)$ is the
Drinfeld double of $G$ (see \seref{s:D(G)}).

First of all, we need to define algebra structure on $\tA$. To do so,
note that it follows from \thref{t:x_g} that for every $g,h\in G$
there exists a unique up to a constant isomorphism of $A$-modules
\begin{equation}\label{e:mut}
\mut_{g,h}\colon X_g\tta X_h\isoto X_{gh}.
\end{equation}

Considering morphisms $X_g\tta X_h\tta X_k\to X_{ghk}$ obtained as
compositions of $\mut$, we get
\begin{equation}\label{e:muassoc}
\mut_{g,hk}\circ(1\tta\mut_{h,k})
=\om(g,h,k) \mut_{gh,k}\circ(\mut_{g,h}\tta 1)
\end{equation}
for some $\om(g,h,k)\in \Ctimes$. One immediately sees that $\om$ is a
3-cocycle on $G$ with values in $\Ctimes$ and that rescaling $\mut$ by
$\mut_{g,h}\mapsto \mut_{g,h}\cdot f(g,h)$ results in replacing $\om$ by a
a cohomological one. Thus, the class $[\om]\in H^3(G, \Ctimes)$ is
well-defined.

To simplify the exposition, in this paper we only consider the
simplest case $\om\equiv 1$. General case is similar but will involve
``twisted'' version of Drinfeld double, as in \cite{DPR}, and will be
discussed elsewhere. Note also that rescaling  $\mut_{g,h}\mapsto
c\mut_{g,h}$ where $c\in \Ctimes$ is independent of $g,h$, does not
change $\om$; thus, without loss of generality we can assume that
$\mut_{1,1}\colon A\tta A\to A$ coincides with the multiplication map
$\mu$.  

\begin{assumption}\label{a:om=1}
  From now on, we assume that $\om\equiv 1$ and $\mut_{1,1}=\mu$. 
\end{assumption}

In this case, the morphism
\begin{equation}
\mut=\bigoplus_{g,h}\mut_{g,h}\colon \tA\tta \tA\to \tA
\end{equation}
is associative. We will also use the the same symbol $\mut$ for the
composition
$$
\tA\ttt \tA\to \tA\tta \tA\xxto{\mut} \tA
$$
where the first morphism is the canonical projection. This morphism
defines on $\tA$ the structure of an associative (but not commutative)
$\C$-algebra.

\begin{lemma}\label{l:murigid}
  Under \asref{a:om=1},
\begin{enumerate}
\item The morphisms $\mut_{1,g}\colon A\ttt X_g\to X_g,
\mut_{g,1}\colon X_g\ttt A\to X_g$ coincide with $\mu_{X_g},
\mu_{X_g}\circ \check R_{A,X_g}^{-1}$ respectively. 
  
\item The morphisms
\begin{equation}
\begin{aligned}
  {}&X_{g^{-1}}\ttt X_g\xxto{\mut}X_1=A\to\one\\
  &\one\injto A=X_1\xxto{(\dim A)\mut^{-1}}X_g\tta X_{g^{-1}}\injto
  X_g\ttt X_{g^{-1}}
\end{aligned} 
\end{equation}
satisfy the rigidity axioms and thus define an isomorphism
\begin{equation}\label{e:x_g*}
X_{g^{-1}}\simeq X_g^*.
\end{equation}
\textup{(}We use the fact that $X\tta Y$ is canonically a direct
summand in $X\ttt Y$, see \cite[Corollary 1.16]{KO}\textup{)}.

\item For all $g$,
  $$
  \frac{\dim X_g}{\dim A}=1.
  $$
\end{enumerate}
\end{lemma}
The proof is left to the reader as an exercise.

Now let us define an action of $G$ on $\tA$. It follows from
\thref{t:x_g} that for every $g,x\in G$ there exists a unique up to a
constant $A$-morphism
\begin{equation}
\ph_{x}(g)\colon X_x^{g^{-1}}\isoto X_{gxg^{-1}}.
\end{equation}
Equivalently, $\ph$ is a $\C$-morphism $X_x\to X_{gxg^{-1}}$ such that
$$
\ph_x(g)\circ \mu\circ (\pi^{-1}_g \ttt 1)= \mu\circ (1\ttt
\ph_x(g))\colon A\ttt X_x\to X_{g^{-1}xg}.
$$
Considering composition $X_x\xxto{\ph_x(g)}
X_{gxg^{-1}}\xxto{\ph_{gxg^{-1}}(h)} X_{hgx(hg)^{-1}}$ and using
uniqueness, we see that
\begin{equation}\label{e:assph}
\ph_{gxg^{-1}}(h)\ph_x(g)=c_x(g,h)\ph_x(hg)
\end{equation}
for some $c_x(g,h)\in \Ctimes$.

In particular, denoting by $Z(x)$ the centralizer of $x$:
$$
Z(x)=\{g\in G\st gx=xg\}
$$
we see that $g\mapsto \ph_x(g)$ defines a projective action of
$Z(x)$ on $X_x$.
\begin{lemma}\label{l:c=1}
  If $\om\equiv 1$, then $\ph_g(x)$ can be chosen so that
  $c_x(g,h)\equiv 1$.
\end{lemma}
\begin{proof}
  
  Define $\ph_x(g)$ by \firef{f:phixg2} (where we used \eqref{e:x_g*}
  to identify $X_g^*\simeq X_{g^{-1}}$). We leave it to the reader to
  check that so defined $\ph$ satisfies the associativity property
  \eqref{e:assph} with $c_x(g,h)\equiv 1$.

\begin{figure}[h]
  $$
  \ph_x(g)=\frac{1}{\dim A}\quad\newfig{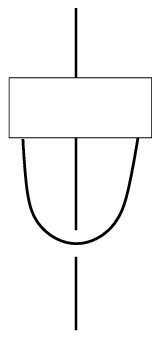}{
    \rput(-0.9, 0.8){$X_x$}
    \rput(-2.05, 1.9){$X_g$}
    \rput(-0.4, 1.9){$X_g^*$}
    \rput(-1.25,2.65){$\mut$}
    }
  $$
\caption{Definition of $\ph_x(g)\colon X_x^{g^{-1}}\isoto X_{gxg^{-1}}$}
\label{f:phixg2}
\end{figure}
\end{proof}

\begin{example}\label{e:phi=pi}
  Let $x=1, X_1=A$. Then it immediately follows from the construction
  given in the proof of \thref{t:twisted1} and \leref{l:murigid} that
  $\ph_1(g)=\frac{1}{\dim X_g}T_g=\pi_g$.
\end{example}

\begin{lemma}\label{l:phimu}
  If $\om\equiv 1$ and $\ph$ is defined as in \leref{l:c=1} then $\ph$
  is compatible with $\mut$, i.e., the following diagram is
  commutative:
\begin{equation}\label{e:phimu}
\begin{CD}
  X_x\ttt X_y @>{\mut}>>X_{xy}\\
  @V{\ph_x(g)\ttt \ph_y(g)}VV                @V{\ph_{xy}(g)}VV\\
  X_{gxg^{-1}}\ttt X_{gyg^{-1}} @>{\mut}>> X_{gxyg^{-1}}
\end{CD}
\end{equation}
\end{lemma}

\begin{remark}
  In general case ($\om\ne 1$), it is easy to see that \eqref{e:phimu}
  is commutative up to a constant factor $\ga_{x,y}(g)\in\Ctimes$. We
  plan to show in a forthcoming paper that both $c_x(g,h)$ and
  $\ga_{x,y}(g)$ can be expressed in terms of $\om$ in a manner
  similar to \cite[Equations 3.5.2, 3.5.3]{DPR}.
\end{remark}

Denote
$$
\ph(g)=\bigoplus_{x}\ph_x(g)\colon \tA\to \tA.
$$
Then we have the following result. (We assume that the reader is
familiar with the definition and properties of the Drinfeld double
$D(G)$ of a finite group $G$; these results are briefly reviewed in
\seref{s:D(G)}.)

\begin{theorem}\label{t:tadg}
  Let $\om\equiv 1$ and $\ph$ defined as in \leref{l:c=1}. Define the
  map $\ph\colon D(G)\to \End_\C\tA$  by
\begin{equation}
\begin{aligned}
  g&\mapsto \ph(g)\\
  \de_h&\mapsto p_h,
\end{aligned}
\end{equation}
where $p_h\colon \bigoplus X_x\to\bigoplus X_x$ is the projection on
$X_h$. Then $\ph$ defines an action of $D(G)$ on $\tA$ by
$\C$-endomorphisms. This action preserves multiplication $\mut$: for
every $x\in D(G)$, $\mut\circ (\ph\ttt \ph)\De(x)=\ph(x)\circ\mut$.
\end{theorem}
\begin{proof}
  Immediately follows from the commutation relations in $D(G)$,
  \leref{l:c=1} and \leref{l:phimu}.
\end{proof}

Note that it follows from \exref{e:phi=pi} that restriction of $\ph$
to $\Cset[G]\subset D(G), A\subset \tA$ coincides with the original
action of $G$ on $A$.

Thus, we have a situation analogous to that of \seref{s:main}: we have
an associative $\C$-algebra $\tA$ on which $D(G)$ acts by
endomorphisms. Analogously to definition in \seref{s:preliminary}, let
$\cd$ be the category with objects: pairs $(M\in\C, \rho\colon
D(G)\to\End_\C(M))$ and morphisms: $\C$-morphisms commuting with the
action of $D(G)$. We have the following results which are parallel to
those given in \seref{s:preliminary} for $\C[G]$.

\begin{lemma}\label{l:cdg}
{\ }\\
\begin{enumerate}
\item $\cd$ has a canonical structure of a rigid monoidal category
  
\item Both $\Rep D(G)$ and $\C$ are naturally subcategories in $\cd$.
  This, in particular, allows us to define the functor of exterior
  tensor product\\
  $\tbox\colon\Rep D(G)\times\cd\to\cd$.
  
\item $\cd$ is a semisimple abelian category with simple objects
  $V_{g,\pi}\tbox L_i$.
  
\item $\cd$ is braided with the commutativity isomorphism
\begin{equation}\label{e:rcd}
\rcd=\check R\circ R^{D(G)}
\end{equation}
where $R^{D(G)}$ is the $R$-matrix of $D(G)$.

\item For $V\in \C, W\in \Rep D[G]$ considered as objects in $\cd$,
  one has $\bigl(\rcd_{VW}\bigr)^2=\id$. 
\end{enumerate}
\end{lemma}

The proof of this lemma is left to the reader as an exercise. Note
that unlike $\C[G]$ case, $\rcd$  does not coincide with the
usual commutativity isomorphism in $\C$.

We can also define the notion of ``$D(G)$ invariants''. Namely, define
functors
\begin{align*}
  \Rep D(G)&\to \vec\\
  V&\mapsto V^{D(G)}=\Hom_{D(G)}(\Cset, V)
\end{align*}
and
\begin{align*}
  \cd&\to \C\\
  \oplus V_i\tbox L_i&\mapsto V_i^{D(G)}\tbox L_i
\end{align*}
(or, more formally, by $\Hom_\C(L, M^{D(G)})=\Hom_{\cd}(L, M)$ where
$L\in\C$ is considered as an object of $\cd$ with trivial action of
$D(G)$). Using semisimplicity of $D(G)$, one easily sees that
$M^{D(G)}$ is canonically a direct summand in $M$, and that for every
$L,M\in \cd$ there is a canonical embedding
$$
L^{D(G)}\ttt M^{D(G)}\injto (L\ttt M)^{D(G)}.
$$

\begin{theorem}
  $\tA$ is an associative commutative algebra in $\cd$ \textup{(}with
  multiplication $\mut$ and action of $D(G)$ defined as in
  \thref{t:tadg}\textup{)}, and $\tA^{D(G)}=A^G=\one$.
\end{theorem} 
\begin{proof}
  The only part which is not obvious is the fact that $\tA$ is
  commutative (with respect to $\rcd$, not $\check R$!), i.e. that the
  composition
  $$
  \tA\ttt \tA\xxto{R^{D(G)}}\tA\ttt \tA\xxto{\check R} \tA\ttt
  \tA\xxto{\mut} \tA
  $$
  coincides with $\mut$. To prove it note that explicit formula
  \eqref{e:rdg} for $R^{D(G)}$ shows that this composition, when
  restricted to $X_{g_1}\ttt X_{g_2}$, is equal to
  $$
  X_{g_1}\ttt X_{g_2} \xxto{1\ttt \ph(g_1)}X_{g_1}\ttt
  X_{g_1g_2g_1^{-1}} \xxto{\check R}X_{g_1g_2g_1^{-1}}\ttt
  X_{g_1}\xxto{\mut} X_{g_1g_2}.
  $$
  Using presentation of $\ph(g_1)$ given in \firef{f:phixg2} and
  associativity of $\mut$, we can rewrite it as shown in
  \firef{f:tacommut}, which shows that it is equal to $\mut$.

\begin{figure}[h]
\begin{align*}
\frac{1}{\dim A}
&  \newfig{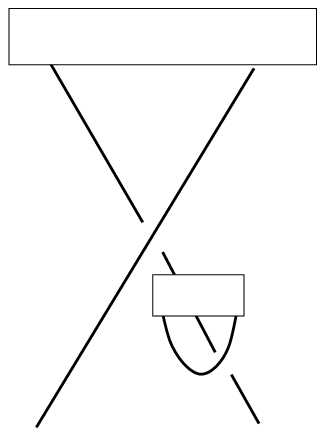}{\rput(-1.6,1.7){$\mut$}\rput(-2,4.35){$\mut$}}
=\frac{1}{\dim A}\newfig{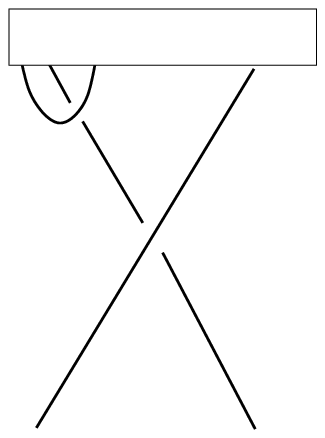}{\rput(-2,4.35){$\mut$}}\\
&=\frac{1}{\dim A}\newfig{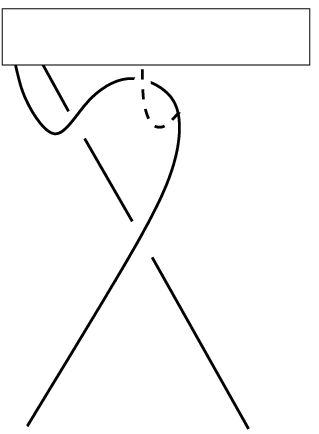}{\rput(-2,4){$\mut$}}
=\newfig{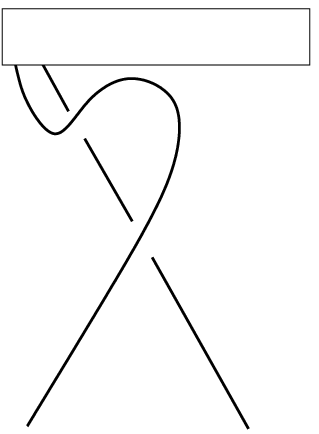}{\rput(-2,4){$\mut$}}
\end{align*}
\caption{Proof of commutativity of $\mut$ with respect to $\rcd$}
\label{f:tacommut}
\end{figure}

\end{proof}

\begin{figure}[h]
\newfig{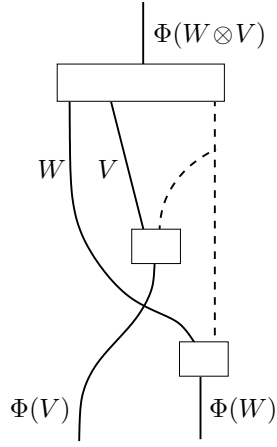}{
\rput(-0.3, 0.8){$\Ph(W)$}
\rput(-3, 0.8){$\Ph(V)$}
\rput(-2.85, 4){$W$}
\rput(-2.1, 4){$V$}
\rput[Bl](-1.5, 5.7){$\Ph(W\sttt V)$}
 }
\caption{Definition of $J$}\label{f:J}
\end{figure}

Now, define the functor $\Phi\colon \Rep D(G)\to \C$ by
\begin{equation}\label{e:tphi}
\Phi(V)=(V\ttt \tA)^{D(G)}
\end{equation}
(cf. with \eqref{e:phi}) and functorial morphism $J\colon \Ph(V)\ttt
\Ph(W)\to \Ph(W\ttt V)$ (note that it reverses the order!) by
\begin{equation}\label{e:tJ}
\begin{aligned}
  {}&(V\ttt \tA)^{D(G)}\ttt (W\ttt \tA)^{D(G)} \isoto 
     \Bigl((V\ttt \tA)^{D(G)}\ttt W\ttt \tA\Bigr)^{D(G)}\\
   &\xxto{\rcd}
     \Bigl(W\ttt(V\ttt \tA)^{D(G)}\ttt \tA\Bigr)^{D(G)}
     \injto (W\ttt V\ttt \tA\ttt \tA)^{D(G)}\\
   &\xxto{\mut}
   (W\ttt V\ttt\tA)^{D(G)}=\Ph(W\ttt V)
\end{aligned}
\end{equation}
(cf. \thref{t:J}). Please note that by \leref{l:cdg}(5), $\rcd$ in the
second line can be replaced by $\bigl(\rcd\bigr)^{-1}$; this will be
used in the future. The definition of $J$ is illustrated in
\firef{f:J}, where --- unlike all previous figures in this paper ---
crossings represent $\rcd$ and not $\check R$, and the dashed line
represents object $\tA$. The boxes are the canonical embeddings
$(V\ttt\tA)^{D(G)}\injto V\ttt \tA$ and projections $V\ttt \tA\surjto
(V\ttt\tA)^{D(G)}$.

\begin{theorem}
  $J$ is compatible with associativity, unit isomorphisms and reverses
  commutativity isomorphism, i.e. $J\circ \check R^{-1} =\Ph(\check
  R^{D(G)})\circ J$.

\end{theorem}
\begin{proof}
  The only one which is not obvious is the commutativity isomorphisms,
  proof of which is shown in \firef{f:jcommut}, which uses the same
  conventions as \firef{f:J}. 
\begin{figure}[h]
\begin{align*}
J\check R&=\newfig{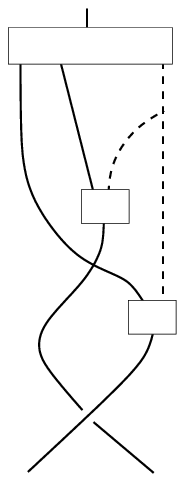}{}\quad
               =\quad\newfig{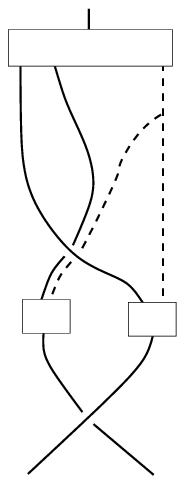}{}\quad
               =\newfig{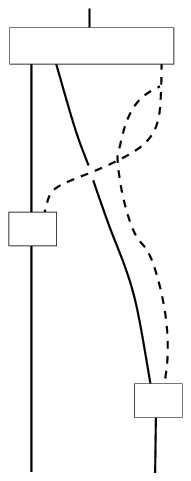}{}\\[10pt]
&=\quad\newfig{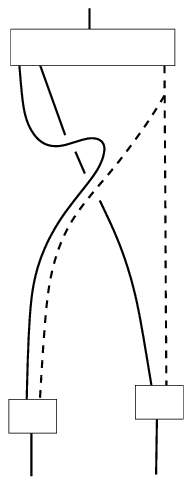}{}\quad
 =\quad\newfig{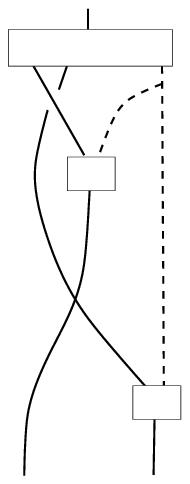}{}=\Ph(\check R^{-1})J.
\end{align*}
\caption{}\label{f:jcommut}
\end{figure}

\end{proof}

Now, repeating with obvious changes the same steps as in
\seref{s:main}, we get the following results. We denote by $\tdg$ the
set of isomorphism classes of irreducible representations of $D(G)$
and for each $\la\in \tdg$ we choose a representative $V_\la$.
\begin{theorem}
\begin{enumerate}
\item For any $\la\in\tdg$, $\Ph(V_\la)$ is either zero or an
  irreducible object in $\C$.
\item If $\la,\mu\in \tdg$ are such that $\Ph(V_\la), \Ph(V_\mu)$ are
  non-zero, and $\la\ne\mu$, then $\Ph(V_\la)\not\simeq V_\mu$.
\item If $\la\in\tdg$ and $\Ph(V_\la)\ne 0$, then $J\colon
  \Ph(V_\la)\ttt \Ph(V)\to \Ph(V\ttt V_\la)$ is an isomorphism.
\item Let $I=\{\la\in\tdg \st \Ph(V_\la)\ne 0\}$, and let $\A$ be the
  abelian subcategory in $D(G)$ generated by $V_\la, \la\in I$. Then
  $\A$ is a subcategory in $\Rep D(G)$ which is closed under taking
  submodules, tensor product, and duality.
\end{enumerate}
\end{theorem}

Combining this last part with \thref{t:main} and \thref{t:rephg}, we
get the following theorem, which is the main result of this paper. 

\begin{theorem}\label{t:main1}
  Assume that action of $G$ is faithful and $\om\equiv 1$. Let
  $H\subset G$ be the normal subgroup defined in \coref{c:H}.  Then
  the functor $\Ph\colon \Rep D(G)\to \C$ defined by \eqref{e:tphi} is
  an equivalence $(\Rep D(G,H))^{op}\isoto \C$.  \textup{(}The
  category $\Rep D(G,H)$ is defined in \eqref{e:dhg} and $(\Rep
  D(G,H))^{op}$ is $\Rep D(G,H)$ with opposite tensor product, see
  \eqref{e:repop}.\textup{)}

  If, in addition, $\C$ is modular then $H=G$ and thus, $\Ph$ is an
  equivalence\\ $(\Rep D(G))^{op}\simeq \C$. 
\end{theorem}

\begin{corollary}\label{c:main2}
Let $\tau\colon \Rep D(G)\to (\Rep D(G))^{op}$ be as defined in
\eqref{e:tau}. Then the composition $\Phi\circ\tau\colon V\mapsto
\Phi(V^\tau)$ is an equivalence of tensor categories $\Rep
D(G,H)\isoto \C$. If $\C$ is modular,  this gives an equivalence
$\Rep D(G)\isoto \C$. 
\end{corollary}

Combining this result with the results of \cite{KO}, we get the
theorem formulated in the introduction.

\begin{corollary}{\ }\\
\begin{enumerate}
\item For every $g\in H$, $X_g$ considered as an object of $\C$ has
  decomposition
  $$
  X_g\simeq \oplus_{\pi}\pi\tbox\Ph(V_{g,\pi}^*)
  $$
  where the sum is over $\pi\in \widehat{Z(g)}$, and $\pi$ is the
  vector space of the representation $\pi$. 
  
\item For every simple object $L_i\in \C$ there exists a unique
  conjugacy class $C$ in $G$ such that $L_i$ appears in decomposition
  of $X_g$ iff $g\in C$.
\end{enumerate}

\end{corollary}
\begin{proof}
Follows from the previous corollary and Equation~\eqref{e:g(xg)}. 
\end{proof}

\section{$D(G)$ revisited}\label{s:D(G)II}
It is instructive to explicitly describe the constructions of the
previous sections, and in particular equivalence $\Rep D(G)\isoto \C$
in the case when $\C=\Rep D(G), A=\F(G)$, so that $G$ acts on $A$ by
automorphisms and $\Rep^0A\isoto \vec$ (see \seref{s:D(G)}). It is
natural to expect that in this case, the functor $\C\isoto \Rep D(G)$
is the identity functor; as we will show, this is indeed so. 

We already have explicit description of the modules $X_g$. The
``multiplication'' map $\mut\colon X_{g_1}\ttt X_{g_2}\to X_{g_1g_2}$
is given by 
\begin{equation*}
e_{x_1, y_1}\ttt e_{x_2, y_2}\mapsto \de_{y_1,y_2}e_{x_1x_2, y_1}.
\end{equation*}
We leave it to the reader to check that this map is indeed a morphism
of $A$-modules, and is associative. Explicit calculation also shows
that the map $\ph_x(g)\colon X_x\to
X_{gxg^{-1}}$ defined by \firef{f:phixg2} can be explicitly written as 
\begin{equation}
\ph_x(g)\colon e_{a,b}\mapsto e_{a,bg^{-1}}
\end{equation}

The object 
$$
\tA=\bigoplus X_g=\bigoplus_{x,y}\Cset e_{x,y}
$$
is a $D(G)\ttt D(G)$-module: the first action is defined as in
\thref{t:tadg} and the second action comes from the fact that each
$X_g$ is an object of $\C= \Rep D(G)$. Explicitly, these two actions
are written as follows:
\begin{alignat}{2}
&\pi^{1}(g)e_{x,y}=e_{x,yg^{-1}},&\quad 
     &\pi^{1}(\de_h) e_{x,y}=\de_{h,y^{-1}xy}e_{x,y},\\
&\pi^{2}(g)e_{x,y}=e_{gxg^{-1},gy},&\quad 
     &\pi^{2}(\de_h) e_{x,y}=\de_{h,x}e_{x,y}.
\end{alignat}

This bi-module admits a more explicit description. Namely, consider
$D(G)$ as a $D(G)\ttt D(G)$-module by $\pi^{1}(x)a=xa,
\pi^{2}(x)a=aS(x)$, where $S$ is the antipode in $D(G)$. Since $D(G)$
is semisimple as an associative algebra, we have
$$
D(G)\isoto\bigoplus V_{g,\pi}\ttt V^*_{g,\pi}
$$
considered as a $D(G)\ttt D(G)$-module: the first copy of $D(G)$
acts on the first factor in the tensor product, the second on the
second. This also shows that the functor $V\mapsto (V\ttt
D(G))^{D(G)}$ (we consider invariants with respect to the action of
$D(G)$ on the tensor product given by $\pi^1$; thus, this space
becomes a $D(G)$-module via $\pi^2$) can be canonically identified
with the identity functor.

How is this related to $\tA$? The answer is given by the following
simple lemma.

\begin{lemma}
  Let $\tau\colon D(G)\to D(G)$ be defined by \eqref{e:tau}, and let
  $D(G)^\tau$ be $D(G)$ considered as a $D(G)\ttt D(G)$-module by
  $\pi^{1}(x)a=xa, \pi^{2}(x)a=aS(\tau x)$. Then the map
\begin{align*}
\tA=\bigoplus\Cset_{x,y}&\to D(G)^\tau\\
                  e_{x,y}&\mapsto y^{-1}\de_x 
\end{align*}
is an isomorphism of $D(G)\ttt D(G)$-modules.
\end{lemma}
The proof is obtained by direct calculation.

Note: the map $\mut\colon \tA\ttt \tA\to \tA$ does {\bf not} coincide
with the multiplication in $D(G)$!

\begin{corollary}
\begin{enumerate}
\item
As a $D(G)\ttt D(G)$-module, $\tA\simeq V^\tau_{g,\pi}\ttt
V^*_{g,\pi}$ 
\item The functor $\Phi\colon V\mapsto (V\ttt \tA)^{D(G)}$ is
  canonically isomorphic to $\tau$. 
\end{enumerate}
\end{corollary}

Thus, the functor $\Phi\circ\tau$ described in \coref{c:main2} can be
identified with the identity functor, as should have been expected. 



\end{document}